\documentclass[11pt]{amsart}
 \usepackage{graphicx, subfigure}
 
\usepackage{amssymb, amsmath, latexsym}
\usepackage{amsfonts}
\allowdisplaybreaks

\catcode`\@=11
\@addtoreset{equation}{section}

\newcommand \RR   {\mathbb{R}}
\newcommand \del  {\partial}

\newcommand \ve   {\varepsilon}

\newcommand \be   {\begin{equation}}
\newcommand \ee   {\end{equation}}

\newtheorem{definition}{Definition}[section]
\newtheorem{theorem}[definition]{Theorem}

\newtheorem{proposition}[definition]{Proposition} 

\newcommand{\bp}{\begin{proof}}
\newcommand{\ep}{\end{proof}}
\newcommand{\fx}{\varphi}
\newcommand{\ft}{\psi}

\begin{document}

\title[The Glimm scheme for hyperbolic systems of balance laws]{A remark on the Glimm scheme for inhomogeneous hyperbolic systems of balance laws}
\author{Cleopatra Christoforou} 
\address{C.C: Department of Mathematics and Statistics, University of Cyprus, 1678 Nicosia, Cyprus}
\urladdr{http://www.mas.ucy.ac.cy/$\sim$kleopatr}
\email{Christoforou.Cleopatra@ucy.ac.cy}


\begin{abstract}
General hyperbolic systems of balance laws with inhomogeneous flux and source are studied. Global existence of entropy weak solutions to the Cauchy problem is established for small $BV$ data under appropriate assumptions on the decay of the flux and the source with respect to space and time. There is neither a hypothesis about equilibrium solution nor about the dependence of the source on the state vector as previous results have assumed.\\
\,\qquad Date: \today 
\end{abstract}

\maketitle
\section{Introduction}\label{S1}

We consider the Cauchy problem to inhomogeneous systems of balance laws in one space dimension
\be\label{S3 balance}
\del_t U+\del_x F(U,x,t)+G(U,x,t)=0
\ee
with initial data
\be\label{S3 in}
U(x,0)=U_0(x).
\ee
Here $x \in \mathbb{R}$, $t\ge 0$ and the state $U(x,t)$ takes values in ${\mathbb{R}}^n$. Also, the flux $F$ and the source $G$ are given smooth functions from 
${\mathbb{R}}^n\times\RR\times\RR_+$ to ${\mathbb{R}}^n$ and for this article, it is important that $F$ and $G$ depend explicitly on $(x,t)$. System \eqref{S3 balance} is strictly hyperbolic,  that is, the Jacobian matrix  $A(U,x,t)=D_UF(U,x,t)$ has $n$ real and  distinct eigenvalues
\be\label{S3 evalues}
 \lambda_1(U,x,t)<\lambda_2(U,x,t)<\ldots<\lambda_n(U,x,t), 
\ee
known as characteristic speeds.

Solutions to homogeneous systems of conservation laws, i.e. $G\equiv 0$ and $F=F(U)$, have been constructed globally in time in the space of bounded variation $BV$ by various methods: the random choice method of Glimm~\cite{Glimm}, the front tracking algorithm~\cite{Bressan} and the vanishing viscosity method~\cite{BB2} under the assumption of small total variation of the initial data. An exposition of the current state of the theory of conservation laws can be found in the manuscripts~\cite{Bressan, Dafermos3, Serre}.

For systems of balance laws with general source term $G$, blow-up of solutions in finite time is expected even when the initial data is of
small variation. In fact,  the presence of the production
term $G$ results to the amplification in time of even small oscillations in the solution.
Because of this feature of hyperbolic balance laws,  one does not expect in general long term stability in BV. Local in
time existence of BV solutions was first established by Dafermos and
Hsiao \cite{DH}, using a modification of the random choice method of Glimm \cite{Glimm}. In addition, in \cite{DH}, Dafermos and Hsiao constructed global BV
solutions under a suitable \emph{dissipativeness assumption} on $G$ using this modified scheme. More precisely, to achieve global existence in~\cite{DH}, a special structure on the dependence of $G$ on the state vector $U$ is assumed as well as existence of a constant equilibrium solution.
Under the same \emph{dissipativeness
assumption} on $G(U,\cdot,\cdot)$ and existence of equilibrium solution, global existence as well as stability results have  been
established via the front tracking method by Amadori and Guerra
\cite{AG, AG1} and via the method of vanishing viscosity by Christoforou \cite{CC06} but only for homogeneous systems, i.e. $F=F(U)$, $G=G(U)$. Recently, this dissipativeness assumption has been relaxed by Dafermos in~\cite{D} to the so--called ``\emph{weak dissipativeness}". For the sake of completeness, we remark that the existence of solutions is also established for other classes of systems of balance laws under different asumptions, for example~\cite{Liu1982, AGG} when $F=F(U)$ and $G=G(U,x)$.

Here, we prove global existence of entropy weak solutions to~\eqref{S3 balance}--\eqref{S3 in} using the modification of Glimm's scheme as constructed in~\cite{DH}, but under different hypotheses.
 It should be emphasized that in this article there are two main characteristics that are not present in the aforemetioned articles: (i) existence of a constant equilibrium solution to~\eqref{S3 balance} is not assumed and (ii) global  existence is established without any conditions on the structure of the source $G$ on $U$ but only through the decay of  inhomogeneity of the flux $F(\cdot,x,t)$ and the source $G(\cdot,x,t)$ with respect to $(x,t)$. Our analysis involves mainly ideas and techniques that are present in the \emph{local existence} result of~\cite{DH}, but implemented under sufficiently rapid decay in the inhomogeneity, we achieve \emph{global} and not only \emph{local} existence. It is important to note that Dafermos and Hsiao in~\cite{DH} comment that global existence can be achieved under appropriate bounds on $F(\cdot,x,t)$ and the source $G(\cdot,x,t)$ and their derivatives. This is also stated in Dafermos's book without proof, cf.~\cite[Chapter $13.9$]{Dafermos3}. Here, this article serves in providing the proof in detail.

It should also be mentioned that another interesting version of the Glimm scheme for inhomogeneous systems of the form~\eqref{S3 balance} has been studied by Hong and LeFloch in~\cite{HL}. The main difference of~\cite{HL} from the aforementioned articles~\cite{DH,AG, D} and this article is in the scheme. More precisely, in~\cite{HL} they solve \emph{generalized Riemann problems} instead of \emph{classical Riemann problems} at every mesh--point of the grid and therefore, the operator splitting is not incorporated in the scheme. Moreover they analyze in detail the new type of nonlinear interactions of the wave patterns, and this leads to a global existence result. Actually, the solution of the \emph{generalized Riemann problems} is carefully approximated and this approximation can be interpreted as the splitting part of the scheme. In this way, one can see similarities between the two versions. Since an approximation is involved, error estimates are crucial in order to ensure the consistency of this version of the method.  Here in this article we use the classical Glimm scheme as constructed by Dafermos and Hsiao~\cite{DH} based on the classical Riemann problem in conjunction with the operator splitting. Last, in both~\cite{HL} and this article, global existence is achieved due to decay of the flux and the source w.r.t. time and space of $L^1$ type.

The structure of this article is as follows: In Section~\ref{S3}, we describe the algorithm of the modified choice method as introduced by Dafermos and Hsiao quoting it from \cite[Section $2$]{DH} for the convenience of the reader. Then, we prove that the approximate sequence of solutions is globally defined and convergences, up to a subsequence, to the entropy weak solution of~\eqref{S3 balance}--\eqref{S3 in}. The main existence result is  Theorem~\ref{thmrandom} stated in Section~\ref{S3.2} and the assumptions are also presented there. The proof of the theorem can be found in Section~\ref{A}. It is technical and involves interaction estimates of elementary waves that yield the compactness of the sequence.

\section{The Modified Random Choice Method}\label{S3}
%
In this section, we describe the algorithm  of the modified random choice method for constructing approximate solutions to balance laws~\eqref{S3 balance}. This algorithm was first introduced in~\cite{DH} and it is based on the random choice method of Glimm~\cite{Glimm} in conjuction with the operator splitting. 

The algorithm is as follows: For strictly hyperbolic system~\eqref{S3 balance} with characteristic speeds~\eqref{S3 evalues}, first construct a staggered grid of mesh points on the upper half plane $(x,t)$ in the following way: Fix a small positive constant $h>0$ corresponding to the space mesh size. Taking $\lambda$ to be a positive constant satisfying
$$\lambda_i(U,x,t)<\lambda,\qquad i=1,\dots n$$
for all $(U,x,t)\in\mathcal{B}\times \RR\times\RR_+$, we choose the time mesh size to be $\Delta t=\lambda^{-1}\,h$. From here and on, $r\in\mathbb{Z}$, $s\in\mathbb{N}_0$ and say $r+s$ is even if $r+s=2k$ for some $k\in\mathbb{Z}$, otherwise it is odd. Next, partition the upper half-plane of the $(x,t)$ plane into strips 
$$\Gamma_s=\{(x,t): x\in\RR,\, s\Delta t\le t<(s+1) \Delta t\},\qquad s=0,1,2,\dots$$
and build the staggered grids of \emph{mesh points} $(x_r, t_s)$, for $r+s$ even, by taking $x_r\doteq rh$ and $t_s\doteq s\Delta t$. As it is described below, the construction of the approximate sequence $U_h(x,t)$ requires solving repeatedly Riemann problems around the mesh points $(x_r, t_s)$ for $r+s$ even over a time interval of length $\Delta t$. Because of the way the time mesh $\Delta t$ is selected,  the waves emanating from the mesh points at a distance $2h$ apart will not interact on a time interval of length $\Delta t$.

In addition, consider a random sequence $\mathcal{A}=\{a_0,a _1,\dots, a_s,\dots\}$ with equidistributed numbers $a_s\in(-1,1)$ and define the \emph{random mesh points} $(y^r_s, t_s)$ with
\be y^r_s=x_r+a_s h,\qquad r+s \text{  odd}. \ee

Assume that $U_h(x,t)$ is known for $x\in\RR$ and $t\in[0,t_s)$ for some positive integer $s$. Define the constant values
\be U^r_s\doteq \lim_{t\to t_s-} U_h(y^r_s -,t),\ee
and  
\be\label{S3.1 hat U}
\hat{U}^r_s\doteq U^r_s-\Delta t\, G(U^r_s, x_r, t_s),
\ee
for every $r$, such that $r+s$ is odd and solve the two equations
\be\label{S3.1:eqs}
F(V^r_s, x_{r+1}, t_s)=F(\hat{U}^r_s, x_r,t_s)=F(W^r_s, x_{r-1}, t_s)\;.
\ee
Now, consider any mesh point $(x_r,t_s)$, i.e. $r+s$ is even, and define $U_h(x,t)$ restricted on the rectangle $\{(x,t): x\in[x_{r-1},x_{r+1}],\, t_s\le t<t_{s+1}\}$ to be the solution to the Riemann problem
\be\label{S3.1:Riem}
\begin{array}{l}
\del_t U_h+\del_x F(U_h, x_r, t_s)=0,\qquad t\ge t_s,\\
U_h(x,t_s)=
\left\{
\begin{array}{ll}
V_s^{r-1} & x<x_r\\
W^{r+1}_s & x>x_r\;.
\end{array}\right.\end{array}
\ee
By solving for every $r$, such that $r+s$ is even, $U_h(x,t)$ is extended for times $t\in(t_s, t_{s+1})$. Assuming that both equations~\eqref{S3.1:eqs} and the Riemann problem~\eqref{S3.1:Riem} are solvable at every point, an approximate solution $U_h(x,t)$ is constructed for $x\in\RR$ and $t>0$. It should be noted that $U_h(x,t)$ admits jump discontinuities across shocks emerging from the mesh points $(x_r, t_s)$ that are part of the Riemann solution as well as across the vertical line segments 
\be
\{ (x,t): x= x_{r-1}, t\in[t_s, t_{s+1})\}
\ee
since $U_h(x_{r-1}-,t)= W^{r-1}_s$ and $U_h(x_{r-1}+,t)=V^{r-1}_s$ for $t\in[t_s, t_{s+1})$ for every $r+s$ even. 
Note that if $F$ is independent of $x$, then this type of discontinuities do not arise. For more details, we refer the reader to~\cite[Sec. $2$]{DH}.

\section{Global Existence}\label{S3.2}
We assume that system~\eqref{S3 balance} is not only stricty hyperbolic, i.e. the Jacobian matrix $D_UF (U,x,t)$ has $n$ real distinct eigenvalues as mentioned in~\eqref{S3 evalues},
which remain uniformly separated, i.e.
\be\label{S3 ass1} |\lambda_i(U,x,t)-\lambda_j(U,x,t)|>A^{-1},\qquad i\ne j,\,\ee
for some positive constant $A$, but also that the characteristic speeds are 
uniformly away  zero, i.e.
\be\label{S3 ass2}  |\lambda_i(U,x,t)|>A^{-1},\qquad i=1,\dots,n,\ee
for all  $U\in\mathcal{B}$, $x\in\RR$ and  $t\ge 0$. Here, $\mathcal{B}$ denotes a ball in $\RR^n$.
By~\eqref{S3 ass2}, we avoid resonance for every $t>0$. Furthermore, we assume that the flux and the source satisfy the uniform bounds
\be\label{S3 ass3} |D_U F(U,x,t)|,\,|D^2_U F(U,x,t)|,\,|D_U G(U,x,t)|\le A\ee
for $U\in\mathcal{B}$, $x\in\RR$ and $t\ge 0$. In addition, there exist a positive function $\fx\in W^{1,1}(\RR)$ with
\be\label{S3 J a} \int_{-\infty}^{\infty} \fx(x)\,dx\le \omega\ee
and a bounded function $\ft\in L^1(0,\infty)$, for which
\be\label{S3 ass4} |F_x (U,x,t)|,\,|G(U,x,t)|\le \omega \ft(t),\ee
\be\label{S3 ass5}|D_U F_x(U,x,t)|,\,|D_U F_t(U,x,t)|,\,|D_U G(U,x,t)\le \omega \ft(t)   ,\ee
\be\label{S3 ass6} |F_x(U,x,t)|, \,|D_U F_x(U,x,t)|,\,|F_{tx}(U,x,t)|,\,|G_x(U,x,t)|\le \fx(x)\,\ft(t), \ee
for every $U\in\mathcal{B}$, $x\in\RR$ and $t\ge 0$.

Under the aforementioned assumptions for system of balance laws~\eqref{S3 balance}, we establish global existence of entropy weak solutions to the Cauchy problem for initial data of small bounded variation using the approximate sequence described in Section~\ref{S3}. This is the main result of this paper:

\begin{theorem}\label{thmrandom}
Under the assumptions~\eqref{S3 ass1}--\eqref{S3 J a} on the flux $F$ and the source $G$ terms, there are positive constants $\omega_0$ and $\delta_0$, depending only on $A$ and the radius $R$ of the ball $\mathcal{B}_R$, such that  
when \eqref{S3 ass4}--\eqref{S3 ass6} are satisfied with $\omega<\omega_0$ and $TV U_0<\delta_0$, then if $U_0\in\mathcal{B}_{R/2}\subset\mathcal{B}_R$, there exists $h_0>0$ small enough for which the family of approximate solutions $\{U_h(x,t)\}_h$, with $h$ being $0<h<h_0$, is globally defined with $(x,t)\in\RR\times[0,\infty)$ for each sequence $\{a_s\}$.
Moreover, there is a converging subsequence $\{U_{h_k}\}$, with $h_k\to0+$ as $k\to\infty$, such that $U_{h_k}$ converges in $L^1_{loc}$ to a function $U$, for almost every sequence $\{a_s\}$, which is the entropy admissible weak solution to the Cauchy problem~\eqref{S3 balance}--\eqref{S3 in}. Furthermore, for each $t>0$, $U(t)$ is a function of bounded variation on $(-\infty,\infty)$ and
\be\label{S3 thmest}
TV_{(-\infty,\infty)}\{U(t)\}\le  C_1 e^{\sigma} (TV\{U_0\}+\omega),
\ee
\be\label{S3 thmest infty}
\|U(t)\|_{L^\infty}\le \|U_0\|_{L^\infty}+C_2 e^{\sigma} (TV\{U_0\}+\omega).
\ee
where $C_1$ and $C_2$ are constants and $\sigma=O(1)\omega\|\psi\|_{L^1[0,\infty)}$.
\end{theorem}

Throughout this article the notation $O(1)$ would mean the term is bounded by a constant.

The consistency of the algorithm follows in the same way as in the local case of~\cite[pp. 476--481]{DH}. Here is stated only the result:
\begin{proposition}[Dafermos and Hsiao~{\cite[pp.476]{DH}}]\label{S3 prop}
 Assume that for each selection of random sequence $\mathcal{A}=\{a_s\}$, the family $\{U_h(x,t)\}_h$, $0<h<h_0$, of approximate solutions is defined for $x\in\RR$ and $t>0$ and has locally bounded variation on $(-\infty,\infty)\times[0,T)$ uniformly in $h$. Then, there is a converging subsequence $\{U_{h_k}\}$, with $h_k\to0+$ as $k\to\infty$, such that, for almost every sequence $\{\alpha_s\}$,
\be
\lim_{k\to\infty} U_{h_k}(x,t)=U(x,t),\qquad \text{a.e.}
\ee
 where $U(x,t)$ is a function of locally bounded variation on $(-\infty,\infty)\times[0,T)$, which is a weak solution to the Cauchy problem~\eqref{S3 balance}--\eqref{S3 in} and satisfies the entropy admissibility criterion.
\end{proposition}

The proof of Theorem~\ref{thmrandom} is technical and follows ideas from the \emph{local existence result} in~\cite{DH}. The proof can be found in Section~\ref{A} and the presentation gives an emphasis to the two main aspects that are not present in previous results~\cite{DH, AG, CC06}, i.e. (i) no equilibrium solution is assumed and (ii) the global existence is due to rapid decay of inhomogeneity and not because of any special structure assumed on $G(U,\cdot,\cdot)$.
We remark that  the aim in the proof is to construct the approximate sequence $U_h$ constructed in Section~\ref{S3}  on the whole upper half plane. The next step is to prove its convergence, up to a subsequence, in $L^1_{loc}$ to the entropy weak solution $U$ of~\eqref{S3 balance},~\eqref{S3 in} for $(x,t)\in\RR\times(0,\infty)$. By Proposition~\ref{S3 prop}, this is accomplished by establishing global bounds of the total variation of $U_h$, \emph{that are uniform in $h$} and time.


\section{Proof of Theorem~\ref{thmrandom}}\label{A}
In this section, we prove uniform bounds on the total variation for the approximate solutions $\{U_h\}$ to~\eqref{S3 balance} as constructed via the modified random choice method in Section~\ref{S3}. Our treatment follows closely the proof of \emph{local} existence in~\cite{DH}, but under the assumptions of Theorem~\ref{thmrandom} on the inhomogeneity of the flux and the source on $(x,t)$ that are stricter than~\cite[Theorem 1.1]{DH}, we obtain a uniform estimate that is \emph{global} in time. For the convenience of the reader, we present all the steps of the proof and refer to~\cite[Section 3]{DH} when necessary. It should be mentioned that this work will require laborious analysis. 

We partition 
$\RR\times[0,\infty)$ 
into a countable set of diamonds $\Delta^r_s$, $r+s$ odd, with vertices at the \emph{random mesh points} $(y^r_s, t_s)$, $( y^{r-1}_{s+1}, t_{s+1})$, $( y^r_{s+2}, t_{s+2})$ and  $( y^{r+1}_{s+1}, t_{s+1})$. The goal is to estimate the strength of the outgoing $\ve$ wave fan which emanates from $(x_r, t_{s+1})$ in terms of the strengths of the incoming wave fans $\alpha$ and $\beta$ emanating from the points $(x_{r-1}, t_s)$ and $(x_{r+1}, t_s)$, respectively.

We consider the wave fan function $\Phi(\tau; \overline{U}, \overline{x}, \overline{t})$, which is constructed by the composition of the $n$ wave fan curves $\Phi_i$ as follows
\be \Phi(\tau; \overline{U}, \overline{x}, \overline{t})=\Phi_n(\tau_n;\Phi_{n-1}(\tau_{n-1};\dots;\Phi_1(\tau_1;\overline{U},\overline{x}, \overline{t}),\dots,\overline{x}, \overline{t}),\overline{x}, \overline{t}),\ee
and is associated with the Riemann solution to the $n\times n$ strictly hyperbolic system of conservation laws
\be\label{A: cons ct} \partial_t U+\partial_x F(U, \overline{x}, \overline{t})=0,\ee
with \emph{homogeneous} flux. It should be noted that $\Phi$ is smooth satisfying
\be \Phi(0,\overline{U}, \overline{x}, \overline{t})=\overline{U},\ee
\be D_\tau\Phi (0,\overline{U}, \overline{x}, \overline{t})=R(\overline{U}, \overline{x}, \overline{t}),\ee
where $R$ is the $n\times n$ matrix of right eigenvectors associated with the Jacobian of the flux $F$. Next, we denote the inverse of $\Phi(\tau; \overline{U}, \overline{x}, \overline{t})$ in  a small neighborhood of $0$ by $\tau=\Omega(V; \overline{U}, \overline{x}, \overline{t})$, which is the Riemann solution to~\eqref{A: cons ct} having as left Riemann data $\overline{U}$ and right Riemann data $V$. The Riemann solution consists of $n+1$ constant states separated by $n$ elementary $i$-waves of strength $\tau_i$.

Examining carefully the diamond $\Delta^r_s$ we see that the following identities hold:
\be\label{A Wrs} W^{r}_s=\Phi(\alpha, U^{r-1}_{s+1}, x_{r-1}, t_s),\ee
\be\label{A Ur1} U^{r+1}_{s+1}= \Phi(\beta, V^r_s, x_{r+1}, t_s),\ee
\be W^{r+1}_{s+1}= \Phi(\ve, V^{r-1}_{s+1}, x_{r}, t_{s+1}).\ee
As in~\cite{DH}, we define $\overline{W}^r_s$, $\overline{U}^{r+1}_{s+1}$ and $\overline{\ve}$ as follows:
 \be\label{A W bar}   \overline{W}^r_s=\Phi(\alpha, U^{r-1}_{s+1}, x_r, t_{s+1}),\ee
 \be\label{A U bar} \overline{U}^{r+1}_{s+1}=\Phi(\beta, \overline{W}^r_s, x_r, t_{s+1}),\ee
   \be \overline{U}^{r+1}_{s+1}=\Phi(\overline{\ve} , U^{r-1}_{s+1}, x_r, t_{s+1}).\ee
By existing results of Glimm~\cite{Glimm} and Liu~\cite{Liu1979}, we deduce
\be\label{A e bar} \overline{\ve}=\alpha+\beta+O(1)\mathcal{D}(\Delta^r_s),\ee
where $\mathcal{D}(\Delta^r_s)$ denotes the total amount of wave interaction in the diamond $\Delta^r_s$, i.e.
$$\mathcal{D}(\Delta^r_s)=\sum_{\text{approaching}} |\alpha^i\beta^j|,
$$
and the summation runs over all pairs of \emph{approaching waves}. For the notion of approaching waves, we refer to the book~\cite[Chapter XIII]{Dafermos3}. Using assumptions~\eqref{S3 ass1}, \eqref{S3 ass2}, \eqref{S3 ass5}, and~\eqref{A Wrs}, \eqref{A Ur1}, we get
\be\label{A Wdiff}|W^r_s-\overline{W}^r_s|\le C h\omega  \,\ft(t_s) |\alpha|.\ee
Here $C$ is a generic constant which may increase as the computation proceeds. The solution of equations~\eqref{S3.1:eqs} yields the existence of an invertible function $\vartheta(U,x,t)$ such that
\be\label{A V} V^r_s=U^r_s-h\vartheta (U^r_s, x_r, t_s)-\Delta t G(U^r_s, x_r, t_s)\ee
\be\label{A W}W^r_s=U^r_s+h\vartheta (U^r_s, x_r, t_s)-\Delta t G(U^r_s, x_r, t_s)\ee
with
\be\label{A theta}\vartheta(U,x,t)=D_UF^{-1}(U,x,t)\,F_x(U,x,t).\ee
By hypotheses~\eqref{S3 ass1}, \eqref{S3 ass2}, \eqref{S3 ass3}, \eqref{S3 ass4}, ~\eqref{S3 ass6}, as well as ~\eqref{A Wdiff}--\eqref{A W}, we obtain the estimate
\be\label{A VW} |V^r_s-\overline{W}^r_s|\le|V^r_s-W^r_s|+|W^r_s-\overline{W}^r_s|\le C\,\omega\, \ft(t_s) h.\ee
Then, by bounds \eqref{S3 ass3}, \eqref{S3 ass5}, equations \eqref{A Ur1},~\eqref{A U bar} and estimate \eqref{A VW}, we get
\be\label{A Uterms} | U^{r+1}_{s+1}- V^r_s-\overline{U}^{r+1}_{s+1}+\overline{W}^r_s|\le C\,\omega\, \ft(t_s) \,h |\beta|.\ee
Following the strategy in~\cite[Section 3]{DH}, we write
\be\label{A W sum} W^{r+1}_{s+1}- V^{r-1}_{s+1}-U^{r+1}_{s+1}+V^r_s-W^r_s+U^{r-1}_{s+1}= h \vartheta^{r+1}_{s+1}+h \vartheta^{r-1}_{s+1}-2h\vartheta^r_s-\Delta t G^{r+1}_{s+1}+\Delta t G^{r-1}_{s+1}\ee
using expressions~\eqref{A V}--\eqref{A W} and then we estimate the right-hand side under hypotheses~\eqref{S3 ass1}--\eqref{S3 ass6}. From now on, we use the abbreviation, $\vartheta^i_j=\vartheta(U^i_j, x_i, t_j)$ and $G^i_j=G(U^i_j, x_i, t_j)$.
First, from~\eqref{A Ur1} and \eqref{A Wrs}, we get
\begin{align}\label{A theta1}
|\vartheta^{r+1}_{s+1}-\vartheta(U^r_s,x_{r+1}, t_{s+1})|&=\left((D_UF^{-1})F_x+D_UF^{-1}\,D_U F_x\right) |U^{r+1}_{s+1}-U^r_s|\nonumber\\
&\le C \fx(x_r) \ft(t_s) h+C\omega\, \ft(t_s)|\beta |+O(h^2),
\end{align}
and, similarly,
\be |\vartheta^{r-1}_{s+1}-\vartheta(U^r_s, x_{r-1}, t_{s+1})|\le C\omega\, \ft(t_s)|\alpha|+Ch\fx(x_r) \ft(t_s) +O(h^2).\ee
Next,
\begin{align} |\vartheta(U^r_s, x_{r+1}, t_{s+1})+\vartheta(U^r_s, x_{r-1}, t_{s+1})&-2\vartheta (U^r_s, x_r, t_s)|=2\,\Delta t\,\vartheta_t(U^r_s, x_{r}, t_{s})+O(h^2)\nonumber\\
&\le C\fx(x_r) \ft(t_s) h+O(h^2),\end{align}
and
\be |G(W^r_s, x_{r-1}, t_{s+1})-G(V^r_s, x_{r+1}, t_{s+1})|\le C \fx(x_r) \ft(t_{s+1}) h+O(h^2).\ee
Last, we have the following estimates
\be
|G(U^{r-1}_{s+1}, x_{r-1}, t_{s+1})-G(W^r_s, x_{r-1}, t_{s+1})|\le D_U G|U^{r-1}_{s+1}-W^r_s|\le\omega\, \ft(t_{s+1})\,|\alpha|,
\ee
\be\label{A G3}
|G(V^{r}_{s}, x_{r+1}, t_{s+1})-G(U^{r+1}_{s+1}, x_{r+1}, t_{s+1})|\le D_U G|V^r_s-U^{r+1}_{s+1}|\le \omega\, \ft(t_{s+1})\,|\beta|,
\ee
which are sharper from those in~\cite{DH}, because of the stricter assumptions~\eqref{S3 ass4}--\eqref{S3 ass6}, which allow $L^1$ decay in time; compare with ($3.27$)--($3.28$) in~\cite[pp.485]{DH}. The above two estimates are important to achieve global existence. Adding~\eqref{A theta1}--\eqref{A G3}, we obtain an estimate of the quantity in~\eqref{A W sum}:
\begin{align}\label{A W sum2}
|W^{r+1}_{s+1}- V^{r-1}_{s+1}-U^{r+1}_{s+1}&+V^r_s-W^r_s+U^{r-1}_{s+1}|\le C h^2\fx(x_r) \ft(t_{s+1})\nonumber\\
&+ C h\omega\, \ft(t_{s+1})(|\alpha|+|\beta|)+O(h^3).
\end{align}
Combining~\eqref{A W sum2} with~\eqref{A Wdiff} and~\eqref{A Uterms}, we arrive at
\begin{align}\label{A 22}
|V^{r-1}_{s+1}-W^{r+1}_{s+1}-U^{r-1}_{s+1}&+\overline{U}^{r+1}_{s+1}|\le  Ch^2 \fx(x_r) \ft(t_{s+1})\nonumber\\
&+C h\omega \ft(t_{s+1})(|\alpha|+|\beta|)+O(h^3).
\end{align}
In addition we find
\be |U^{r-1}_{s+1}-V^{r-1}_{s+1}|=|h\vartheta^{r-1}_{s+1}+\Delta t G^{r-1}_{s+1}|\le C h \omega\ft(t_{s+1}),\ee
\be\label{A 24}|U^{r-1}_{s+1}-\overline{U}^{r+1}_{s+1}|\le|U^{r-1}_{s+1}-\overline{W}^r_s|+|\overline{W}^r_s-\overline{U}^{r+1}_{s+1}| \le  C(|\alpha|+|\beta|),\ee
where we have used~\eqref{A W bar} --\eqref{A U bar}.

To get the approximate balance laws of the elementary waves, we introduce the map
\be\Psi(W; U,x,t)=\Omega(U+W; U,x,t)\ee
defined for $(U,x,t)\in\mathcal{B}\times\RR\times[0,\infty)$ and $W\in\RR^n$ with $|W|$ small. It is easy to verify the properties
$\Psi(0; U,x,t)=0$, $|D_U\Psi(W; U,x,t)|\le C|W|$ and $|D_W\Psi(W; U,x,t)|\le C$. 
Then we estimate the difference
\begin{align} |\ve-\overline{\ve}|=|\Psi&(W^{r+1}_{s+1}-V^{r-1}_{s+1}  ; V^{r-1}_{s+1}, x_r, t_{s+1})-\Psi (\overline{U}^{r+1}_{s+1}-U^{r-1}_{s+1} ; U^{r-1}_{s+1}, x_r, t_{s+1})|\nonumber\\
&\le Ch\omega \ft(t_{s+1}) (|\alpha|+|\beta|)+ Ch^2 \fx(x_r) \ft(t_{s+1})+O(h^3),\end{align}
using~\eqref{A 22}--\eqref{A 24}. Hence, by~\eqref{A e bar}, we deduce
\begin{align}\label{A apprx bal} |\ve-(\alpha+\beta)|\le& Ch\omega \ft(t_{s+1}) (|\alpha|+|\beta|)+ Ch^2 \fx(x_r) \ft(t_{s+1})\nonumber\\
&+O(h^3)+C_0 \mathcal{D}(\Delta^r_s).\end{align}
In the following, we use the important estimate~\eqref{A apprx bal} to prove bounds on the total variation of $U_h$.

As in~\cite{Glimm}, we employ the so-called  mesh curves $I$, the associated with them functionals 
\be\label{A L} \mathcal{L}(I)=\sum\{ |\alpha|: {\alpha \text{ crosses } I}\}\ee
\be \mathcal{Q}(I)=\sum\{|\alpha||\beta|: {\alpha, \beta \text{ approaching and  cross } I}\}\ee
and the \emph{Glimm functional}
\be\label{A G} \mathcal{G}(I)=\mathcal{L}(I)+2 C_0 \mathcal{Q}(I).\ee
with $C_0>0$ the constant in~\eqref{A apprx bal}. Using the partial ordering of the set of mesh curves $I$ used in this context, we consider an immediate $(r,s)$ successor $J'$ to $J$, then using~\eqref{A apprx bal}, we have
\begin{align}  \mathcal{L}(J')&\le \mathcal{L}(J)+[Ch\omega(|\alpha|+|\beta|)+C h^2 \fx(x_r)]\ft(t_{s+1})\nonumber\\
&+O(h^3)+ C_0 \mathcal{D}(\Delta^r_s),
\end{align}
\begin{align} \mathcal{Q}(J')&\le  \mathcal{Q}(J)+ \mathcal{L}(J)\big[ [Ch(|\alpha|+|\beta|)+Ch^2 \fx(x_r)]\ft(t_{s+1})\nonumber\\
&+O(h^3)+C_0 \mathcal{D}(\Delta^r_s)\big]-\mathcal{D}(\Delta^r_s),
\end{align}
and
\begin{align}\label{A G bound}
 \mathcal{G}(J')&\le  \mathcal{G}(J)-C_0 [1-2C_0 \mathcal{L}(J)]\mathcal{D}(\Delta^r_s)\nonumber\\
& +[1+2 C_0 \mathcal{L}(J)] \left[ [Ch\omega(|\alpha|+|\beta|)+ Ch^2 \fx(x_r)]\ft(t_{s+1})+O(h^3)\right].
\end{align}

Now, we follow the standard argument in Glimm's scheme~\cite{Glimm}: Let $J_s$ be a mesh curve which originates and terminates on the $x$--axis confined in the strip $t\in[0,t_s]$. Suppose that $\mathcal{L}(J)\le (2C_0)^{-1}$ for every $J\le J_s$. Then, consider a decreasing sequence of mesh curves $\{J_k\}$, $k=0,1,\dots, k$, $J_s\ge J_{s-1}\ge\dots\ge J_1$, such that the mesh curves $J_k$ and $J_{k+1}$ share all nodes in the strip $[0,t_k]$ and if $J_{k+1}$ contains the node $(y^r_{k+1}, t_{k+1})$, then $J_k$ contains the node $(y^r_{k-1}, t_{k-1})$. By virtue of~\eqref{A G bound} and the integrability of $\varphi$, i.e.~\eqref{S3 J a},
\be
 \mathcal{G}(J_{s+1})\le \mathcal{G}(J_{s})+[Ch\omega  \mathcal{G}(J_{s})+Ch \omega]\ft(t_{s+1})+O(h^2).
\ee
Hence we have 
\be  \mathcal{G}(J_{s})\le  e^{\sigma}\mathcal{G}(J_{1})+\omega (e^{\sigma}-1) +O(h)\ee
with $\sigma=O(1)\omega\|\psi\|_{L^1[0,\infty)}$.
By account of~\eqref{A L}--\eqref{A G}, it follows $\mathcal{G}(J_{1})\le C_1 TV\{U_0\}$. Then as in~\cite[pp.487]{DH}, we get
\be\label{A TV end} TV_{J_s} \{U_h\}\le C e^{\sigma} (TV\{U_0\}+\omega)+O(h),\ee
and observe
\be \sup_{J_s} |U_h|\le \sup|U_0(x)|+TV_{J_s} \{U_h\}.\ee
Note that the upper bound of~\eqref{A TV end} is independent of $s$ and compare with bound ($3.47$) in~\cite[pp.487]{DH}. By choosing the parameters $\omega_0$ and $\delta_0$ sufficiently small we can extend the approximate solution $U_h$ for all times. For more detials on this argument, we refer to~\cite[pp.487--488]{DH}. By Proposition~\ref{S3 prop}, the proof of Theorem~\ref{thmrandom} is complete.

\section*{Acknowledgements}     
Christoforou was partially supported by the Start-Up fund 2011-2013 from University of Cyprus. 



\begin{thebibliography}{10}
\bibitem{AGG} D. Amadori, L. Gosse and G. Guerra,  Global BV entropy solutions and uniqueness for hyperbolic systems of balance laws, \emph{Arch. Ration. Mech. Anal.} {\bf 162} (4) (2002) 327--366.

\bibitem{AG} D. Amadori and G. Guerra, 
Global weak solutions for systems of balance laws,
{\em Appl. Math. Lett.}, {\bf 12} (6) (1999) 123--127.

\bibitem{AG1} D. Amadori and G. Guerra, 
Uniqueness and continuous dependence for systems of balance laws with dissipation,
{\em Nonlinear Anal.}, {\bf 49} (2002) 987--1014.


 \bibitem{BB2} S. Bianchini and A. Bressan, 
Vanishing viscosity solutions of nonlinear hyperbolic systems,
{\em Ann. of Math.} {\bf 161} (2005) 223--342.  

\bibitem{Bressan} A. Bressan, 
{\sl Hyperbolic Systems of Conservation Laws: The One-Dimensional Cauchy Problem,} 
Oxford Univ. Press, 2000.


\bibitem{CC06} C. Christoforou,  Hyperbolic systems of balance laws via vanishing viscosity,
{\em Journal of Differential Equations} {\bf 221} (2006), (2), 470--541.

\bibitem{CC06b} C. Christoforou, Uniqueness and sharp estimates on solutions to hyperbolic systems with dissipative source, {\em Commun. Partial Differential Equations,} {\bf 31}, (2006), (12), 1825--1839. 


\bibitem{Dafermos3} C.~M. Dafermos,
{\em Hyperbolic Conservation Laws in Continuum Physics,}  
Grundlehren Math. Wissenschaften Series {\bf 325}, Springer-Verlag, 2010.

\bibitem{D} C.~M. Dafermos,
Hyperbolic systems of balance laws with weak dissipation, {\em J. Hyp. Diff. Equations} {\bf 3} (2006) (3), 505--527.

\bibitem{DH} C.~M. Dafermos and L. Hsiao, 
Hyperbolic systems of balance laws with inhomogeneity and dissipation, 
{\em Indiana U. Math. J.}, {\bf 31} (1982) 471--491.


\bibitem{Glimm} J. Glimm, Solutions in the large for nonlinear hyperbolic systems of equations, 
{\em Comm. Pure Appl. Math.} {\bf 4} (1965), 697--715.

\bibitem{HL} J.~Hong and P.G.~LeFloch, A version of the Glimm method based on generalized Riemann problems, \emph{Port. Math.} {\bf 64} (2007), 199--236.  

\bibitem{Liu1979} T.-P.\ Liu, Quasilinear hyperbolic systems. {\em Comm. \ Math.\ Phys.\ }
{\bf 68}, {\bf  no. 2} (1979),  142--172. 

\bibitem{Liu1982} T.-P.\ Liu, Nonlinear stability and instability of transonic flows through a nozzle. {\em Comm. \ Math.\ Phys.\ }
{\bf 83}, {\bf  no. 2} (1982),  243-260. 


\bibitem{Serre} D. Serre, 
{\em Systems of Conservation Laws, I,II}, 
Cambridge University Press, Cambridge, 1999.





\end{thebibliography}
\end{document}